\documentclass[12pt,letterpaper,twoside]{article}
\usepackage[utf8]{inputenc}  

\usepackage[top=1in,bottom=1in,left=1in,right=1in,includehead]{geometry}
\usepackage[hyphens]{url}
\usepackage{calc}
\usepackage{lastpage}

\usepackage[sc]{mathpazo}
\usepackage[T1]{fontenc}
\usepackage{amsmath,amssymb,amsthm}
\numberwithin{equation}{section}
\usepackage{graphicx}
\usepackage{caption}
\usepackage{enumitem}

\usepackage[dvipsnames,table]{xcolor}
\definecolor{rub-blau}{RGB}{0,53,96}
\definecolor{rub-grun}{RGB}{141, 174, 16}
\definecolor{fast-rub-blau}{RGB}{160,173,189}

\usepackage[numbers,sort&compress]{natbib}

\usepackage[colorlinks=true,citecolor=fast-rub-blau,linkcolor=rub-blau,urlcolor=rub-grun,breaklinks=true]{hyperref}

\usepackage{doi}
\usepackage{url}

\makeatletter

\bibpunct{[}{]}{,}{n}{,}{,}

\newcommand{\biburl}[1]{\href{#1}{Link}}

\makeatother

\usepackage[capitalize,nameinlink,noabbrev]{cleveref}

\usepackage[titletoc,title,page]{appendix}

\newcounter{example}

\usepackage{tcolorbox}
\usepackage{float}
\usepackage{booktabs}
\usepackage{tabularx}
\usepackage{lipsum}
\usepackage{microtype}

\usepackage{listings}

\lstdefinestyle{pythonstyle}{
    language=Python,
    basicstyle=\scriptsize\ttfamily,
    keywordstyle=\color{blue},
    commentstyle=\color{gray},
    stringstyle=\color{gray},
    numberstyle=\tiny\color{gray},
    breaklines=true,
    breakatwhitespace=false,
    tabsize=2,
    showspaces=false,
    showstringspaces=false,
    showtabs=false,
    frame=none,
    rulecolor=\color{black},
   backgroundcolor=\color{gray!5},
}

\lstdefinestyle{markdownstyle}{
    basicstyle=\scriptsize\ttfamily,
    breaklines=true,
    frame=none,
    backgroundcolor=\color{gray!5},
    tabsize=2,
    showspaces=false,
    showstringspaces=false,
    showtabs=false,
    escapeinside={||},  
}

\newtcolorbox{markdownbox}[1][Markdown Code]{
    colback=gray!5,
    colframe=gray!70!black,
    fonttitle=\bfseries,
    rounded corners,
    boxrule=1pt,
    left=0pt,
    right=0pt,
    top=1pt,
    bottom=1pt
}

\newtcolorbox{pythonbox}{
    colback=gray!5,
    colframe=gray!70!black,
    rounded corners,
    boxrule=1pt,
    left=-1pt,
    right=1pt,
    top=0pt,
    bottom=-1pt
}

\renewenvironment{abstract}{%
  \list{}{%
    \leftmargin 25pt
    \rightmargin 25pt
  }\item[]
  \small
  \noindent\textbf{Abstract.} 
}{\endlist}

\newcommand{\keywords}[1]{%
  \par\vspace{6pt}
  \small\noindent\textbf{Keywords:} #1
}

\usepackage{fancyhdr}
\fancypagestyle{firstpage}{%
  \fancyhf{}
  \fancyhead[C]{\footnotesize\textit{Your Conference Proceedings Info Here}}
}

\makeatletter
\newcommand\addressmark[1]{\textsuperscript{#1}}

\def\@title{\relax}
\def\@titlehead{\relax}
\def\title{\@ifnextchar[{\@gettitle@two}{\@gettitle@one}}
\def\@gettitle@one#1{\@gettitle@two[#1]{#1}}
\def\@gettitle@two[#1]#2{\gdef\@title{#2}\gdef\@titlehead{#1}}

\def\@author{\relax}
\def\@authorhead{\relax}
\gdef\@authorand{and}
\gdef\and{\ignorespaces \space \@authorand \space \ignorespaces}
\def\author{\@ifnextchar[{\@getauthor@two}{\@getauthor@one}}
\def\@getauthor@one#1{\@getauthor@two[#1]{#1}}
\def\@getauthor@two[#1]#2{\gdef\@author{#2}\gdef\@authorhead{#1}}

\def\@address{\relax}
\def\address#1{\gdef\@address{#1}}

\def\@abstract{\relax}
\long\def\abstract#1{\long\gdef\@abstract{#1}}
\def\@keywords{\relax}
\renewcommand\keywords[1]{\gdef\@keywords{#1}}  

\renewcommand{\maketitle}{%
  \begingroup
  \thispagestyle{empty}%
  \renewcommand\thefootnote{\fnsymbol{footnote}}%
  \begin{center}%
    {\LARGE\hyphenpenalty=10000 \@title \par}%
    \vskip 12pt%
    {\fontsize{13.5pt}{16pt}\selectfont \@author \par}%
    \end{center}%
    \vskip 6pt%
    {\small\itshape \@address \par}%
  
  \@thanks
  \vskip 12pt%
  \ifx\@abstract\relax\else
    {\small\noindent\textbf{Abstract.} \@abstract\par}
    \vspace{6pt}%
  \fi
  \ifx\@keywords\relax\else
    {\small\noindent\textbf{Keywords:} \@keywords\par}
  \fi
  \setcounter{footnote}{0}%
  \renewcommand\thefootnote{\arabic{footnote}}%
  \endgroup
}
\makeatother

\title[Even with AI, Bijection Discovery is Still Hard]{Even with AI, Bijection Discovery is Still Hard: The Opportunities and Challenges of OpenEvolve for Novel Bijection Construction}

\author[D. Brown, J. He, H. Jenne, H. Kvinge, and M. Vargas]{Davis Brown\addressmark{1,3}, Jesse He\addressmark{1,2}, Helen Jenne\thanks{\href{mailto:helen.jenne@pnnl.gov}{helen.jenne@pnnl.gov}}\addressmark{1}, Henry Kvinge\addressmark{1,4}, \and Max Vargas\addressmark{1}}

\address{\addressmark{1} Pacific Northwest National Laboratory, \addressmark{2} University of California San Diego, \\\addressmark{3} University of Pennsylvania, \addressmark{4} University of Washington}

\abstract{Evolutionary program synthesis systems such as AlphaEvolve, OpenEvolve, and ShinkaEvolve offer a new approach to AI-assisted mathematical discovery. These systems utilize teams of large language models (LLMs) to generate candidate solutions to a problem as human readable code. These candidate solutions are then `evolved' with the goal of improving them beyond what an LLM can produce in a single shot. While existing mathematical applications have mostly focused on problems of establishing bounds (e.g., sphere packing), the program synthesis approach is well suited to any problem where the solution takes the form of an explicit construction. With this in mind, in this paper we explore the use of OpenEvolve for combinatorial bijection discovery. We describe the results of applying OpenEvolve to three bijection construction problems involving Dyck paths, two of which are known and one of which is open. We find that while systems like OpenEvolve show promise as a valuable tool for combinatorialists, the problem of finding novel, research-level bijections remains a challenging task for current frontier systems, reinforcing the need for human mathematicians in the loop. We describe some lessons learned for others in the field interested in exploring the use of these systems.}


\begin{document}

\maketitle

\section{Introduction}
\label{sect:intro}

Evolutionary systems have recently emerged as powerful tools for scientific discovery, leveraging the ability of large language models (LLMs) to generate and modify code. For example, FunSearch was successfully applied to problems in combinatorial optimization such as the cap set problem \cite{ellenberg2025generative, romera2024mathematical}. 
Its successor AlphaEvolve  discovered faster matrix multiplication algorithms, improvements to the minimum overlap problem, and an improved construction of kissing numbers in dimension $11$ \cite{alphaevolve}. These successes were extended in \cite{georgiev2025mathematical} to include advances in finite field Kakeya and Nikodym sets, Crouzeix’s conjecture, and the moving soft problem. These and other evolutionary systems \cite{cheng2025barbarians,lange2025shinkaevolve, liu2024llm4ad} have also been applied to problems in systems research and algorithmic design. 

Evolutionary program synthesis approaches use LLMs to mutate and recombine previously generated code. Using code as a medium to develop machine learning-driven solutions facilitates interpretability and may bias towards simpler solutions. Current systems of this kind can be applied to problems satisfying two requirements: (1) the problem can be formulated so that candidate solutions are expressible in a common programming language (e.g., Python), and (2) candidates can be verified automatically.  

Compared to other areas of math, problems in algebraic combinatorics are particularly well suited for such systems. The objects of interest are often straightforward to represent on a computer, and the software library Sage \cite{sage} simplifies their manipulation and analysis. Furthermore, the likely presence of Sage in LLM training data means models may be able to leverage this code. \textbf{In this paper, we focus on the problem of {\em bijection discovery}: given finite sets $A$ and $B$ (often parameterized by a positive integer $n$), the goal is to produce a mathematically meaningful Python function defining an algorithmic bijection between $A$ and $B$.} We choose problems where the sets $A$ and $B$ are easily generated in Python to enable automatic verification of the proposed bijections.

There is reason for optimism that evolutionary systems could successfully discover bijections. LLM training sets likely contain extensive descriptions of key combinatorial bijections in both natural language (via arXiv and MathOverflow) and code (via Sage), providing prototypes for the kinds of constructions that mathematicians value. Furthermore, LLMs can iterate through potential solutions far more rapidly than humans and draw from a vast breadth of mathematical knowledge. However, whether they are sufficiently creative for this type of discovery remains an open question.

In this paper we present case studies applying the system OpenEvolve \cite{openevolve} to three bijection discovery problems: two with known solutions and one open problem. All problems involve Dyck paths, chosen because their numerous bijections with other objects counted by the Catalan numbers are well-documented in both literature and code, ensuring strong representation in model training data. The problems are:

\begin{itemize}
\item[(1)] {\bf  Bijection between odd-diagonal avoiding paths and Dyck paths}: North East lattice paths from $(0, 0)$ to $(2n, 2n)$ that avoid the points $(2i-1, 2i-1)$, $1 \leq i \leq n$ are enumerated by the Catalan number $C_{2n}$ \cite[Problem A4] {stanley2015catalan}. Stanley \cite{stanley2015catalan} rates this problem as difficulty 2+, which he describes as ``about the hardest problem that could be reasonably assigned to a class of graduate students''. While a bijective proof is known, direct LLM prompting failed to produce the solution, making it an appropriate benchmark for evolutionary methods.

\item[(2)] {\bf Bijection between 321-avoiding permutations and Dyck paths}: Although there are at least three known bijections between Dyck paths
to 321-avoiding permutations \cite{callan2007bijections}, they are nontrivial. 
We selected this problem based on preliminary experiments where we prompted LLMs to produce bijections between Catalan objects without evolutionary refinement. While even weaker models successfully generated some bijections (such as between Dyck paths and 213-avoiding permutations), only GPT-5 produced a bijection with 321-avoiding permutations (the Billey--Jockusch--Stanley bijection \cite{billey1993some}). The difficulty of this problem suggests it could be a useful case study for OpenEvolve when used with other LLMs. 

\item[(3)] {\bf Area-bounce exchanging bijection on Dyck paths}: A recent paper \cite{ayyer2025area} defines an area-bounce exchanging bijection on a subset of Dyck paths. This problem is well-suited for OpenEvolve because the existing implementation \cite{qtcatalan-bijection} serves as a starting point, and extending this bijection to a broader class of Dyck paths would be progress on an open research question. 

\end{itemize}

Our experiments yield mixed but instructive results.
OpenEvolve successfully discovers the known bijection for problem (1), but surprisingly
fails to find a bijection for problem (2) under a variety of different configurations (despite GPT-5 finding a bijection when prompted directly). Similarly, for problem (3), the system does not find a solution. These failures provide guidance for future applications and inform this paper's focus on synthesizing lessons for mathematicians considering using evolutionary approaches in combinatorics research. The paper is organized as follows: Section 2 describes the OpenEvolve framework, Section 3 presents our case studies, and Section 4 discusses lessons learned for problem selection, objective function design, and realistic expectations when using these tools in mathematical research.

\section{Background}
\label{sec:background}
\subsection{Evolutionary program synthesis}

Evolutionary program synthesis systems like OpenEvolve use LLMs to generate computer programs to optimize a specific objective function. For example, if the objective is sphere packing, each generated program might specify a configuration of spheres which can each be judged by their sphere density. A population of these program is generated and then evolves as sub-optimal solutions are discarded, solutions are mutated (also via LLMs) to introduce novel features, and new programs are generated. Given the power of modern machine learning, it is reasonable to ask why one should go through the intermediary of code rather than optimizing for a solution directly. Forcing AI systems to produce solutions as computer programs leads to intrinsic interpretability since an expert can directly examine the code, a critical feature when mathematical discovery is the goal. Being able to run a solution as a computer program also allows for verification of correctness (though not at the level of a proof), thus avoiding the common situation where LLMs produce extremely plausible sounding natural language solutions that contain subtle errors. Moreover, requiring solutions to be written in code injects a powerful simplicity bias which avoids `bag of heuristic' solutions where a model simply enumerates an exhaustive list of conditions that offer no mathematical insight \cite{nikankin2024arithmetic}.

The basic OpenEvolve algorithm proceeds as follows: 
\begin{itemize}
    \item[(1)] Start with an initial Python program that takes as input an elements of a set $A(n)$ for $n \in \mathbb{Z}_{>0}$, and outputs an object in set $B(n)$. 
    \item[(2)] Evaluate this program using an objective function designed to measure how close it is to being a bijection. 
    \item[(3)] Add the program and evaluation to a database (`\emph{population}') of programs. 
    \item[(4)] Prompt an LLM to generate modifications to the program aimed at improving its evaluation. 
    \item[(5)] Apply the LLM's modifications to create a new program and evaluate it.  
    \item[(6)] Repeat steps (4) and (5) for a user-specified number of iterations so as to \emph{evolve} new solutions. 
\end{itemize}

\paragraph{The initial Python program.} For each problem, we started with multiple initial Python programs, evolved in parallel using independent populations (called `islands'). We obtained the initial Python programs by prompting an LLM with a description of the problem, sampling with the default hyperparameters in the chat interface.

\paragraph{Program evaluation.} We evaluated programs by (1) running the program on all objects in set $A(n)$ and measuring statistics that quantified how closely it approximated a bijection and (2) prompting an LLM to evaluate the program.


{\em Empirical evaluation for a fixed $n$.} For a small $n$, we scored each program $f$ by its proximity to being a bijection. To quantify nearness of a map to being a surjection we calculate the \emph{surjectivity score} which we define as the ratio $|f(A(n))|/|B(n)|$. To calculate nearness of a map to be being injective we calculate the \emph{injectivity score}, which we define as the ratio of unique elements in $f(A(n))$ to the total elements in $f(A(n))$. Finally, to penalize $f$ that mapped to the wrong domain we also calculated the fraction of $f(A(n))$ that actually belong to $B(n)$. We call this the \emph{validity score}. We averaged the three scores to obtain a single combined score.

{\em LLM evaluation.} While we use the word ``bijections'' in this work and in the broader community, it is understood that what we actually mean is something more specific. Assigning a random pairing between elements in two sets is of limited value. What we really want is a map $A(n) \rightarrow B(n)$ that both generalizes to most $n$ and also yields insight into $A(n)$ and $B(n)$. To get at these less quantifiable properties, we used an LLM judge to perform additional checks:
\begin{itemize}[noitemsep]
\item[(1)]Check for ``cheating''. We observed that LLMs readily exploit shortcuts to achieve high scores. For example, models sometimes generated all objects from sets $A(n)$ and $B(n)$ and established a one-to-one correspondence using an index-based mapping. The LLM evaluator was instructed to give scores of 0.0 in these cases. 
\item[(2)] Evaluate whether the code implemented the algorithm described in the program's docstring\footnote{When generating programs, LLMs must produce a docstring describing what the program does.}. The purpose of this evaluation was to differentiate programs that scored low in the empirical evaluation due to a flawed algorithm, and programs that scored low due to implementation errors. The LLM was instructed to give a score on a scale of 0.0 to 1.0 along with an argument that the code matches the description in the docstring, or a counterexample. 
\item[(3)] To impose a simplicity prior (on a scale of 0.0 to 1.0), the LLM was prompted to give high scores to programs that are simple, intuitive, and reveal a structural correspondence between the two paths, and low scores to programs that seemed convoluted, arbitrary, or involved a lot of case work. 
\end{itemize}

The LLM was also prompted to provide reasoning for its scores and suggestions for improvement. These responses were used in the prompts in the evolution step, in order to give more guidance than a numeric score. 
The program's final score was a weighted average of the averaged LLM scores and the empirical evaluation.

\paragraph{Prompt for revision.} In the evolution step,
an LLM is given a prompt that explains the desired bijection (including details around $A(n)$ and $B(n)$), provides hints, describes the evaluation criteria, specifies the program to edited, and highlights other high-scoring and/or diverse programs in the current database to take inspiration from.


\paragraph{Program sampling} OpenEvolve uses the Multi-Dimensional Archive of Phenotypic Elites (MAP-Elites) algorithm \cite{Mouret2015IlluminatingSS} to sample programs. The purpose of applying MAP-Elites in this context is to maintain a diverse population of programs to sample from, reducing the likelihood that the system will get stuck in local minima. The algorithm works by discretizing user-defined features, (such as code complexity or difficulty measures) to create a grid of cells, with each cell containing at most one program. When a new program is produced, its corresponding grid cell is determined, and it replaces 
any existing program in that cell if the new one achieves a higher combined score from the empirical and LLM evaluations. We select a program to evolve by randomly sampling from the MAP-Elites grid, but more sophisticated sampling strategies are also possible.

\subsection{Related Work}

A number of recent works have applied evolutionary program synthesis tools to problems in research-level mathematics. In a follow-on to \cite{alphaevolve}, \cite{georgiev2025mathematical} applied AlphaEvolve (the closed-source inspiration for OpenEvolve) to 67 different problems in mathematics with impressive success: in many cases the system was able to discover the best known result and sometimes improved upon it. Notably however, most problems in this work amount to establishing tighter upper or lower bounds on a numerical quantity. Several of their observations are applicable to our use-case and will be discussed further in Section~\ref{sec:lessons}. 

There are other computational tools to assist mathematicians with bijection discovery, including databases like the OEIS \cite{oeis} and FindStat \cite{findstat}, and more recently the Bijectionist's Toolkit \cite{grosz2023bijectionist}. These tools have been applied to study homomesy \cite{dowling2023homomesy, elder2024homomesies} and cyclic sieving \cite{adams2024cyclic}.
We see these approaches as complementary and envision a future iteration where LLMs can use these tools in their search for bijections.

\section{Case studies}

\subsection{Odd-diagonal avoiding paths and Dyck paths}

\paragraph{Problem.}
Produce a program implementing a bijective map from NE lattice paths from $(0, 0)$ to $(2n, 2n)$ that avoid the points
$(2i-1, 2i-1)$ to Dyck paths from $(0, 0)$ to $(2n, 2n)$. 

\begin{figure}[!htb]
    \centering
\begin{minipage}{.475\textwidth}
\centering
\includegraphics[width=.95\textwidth]{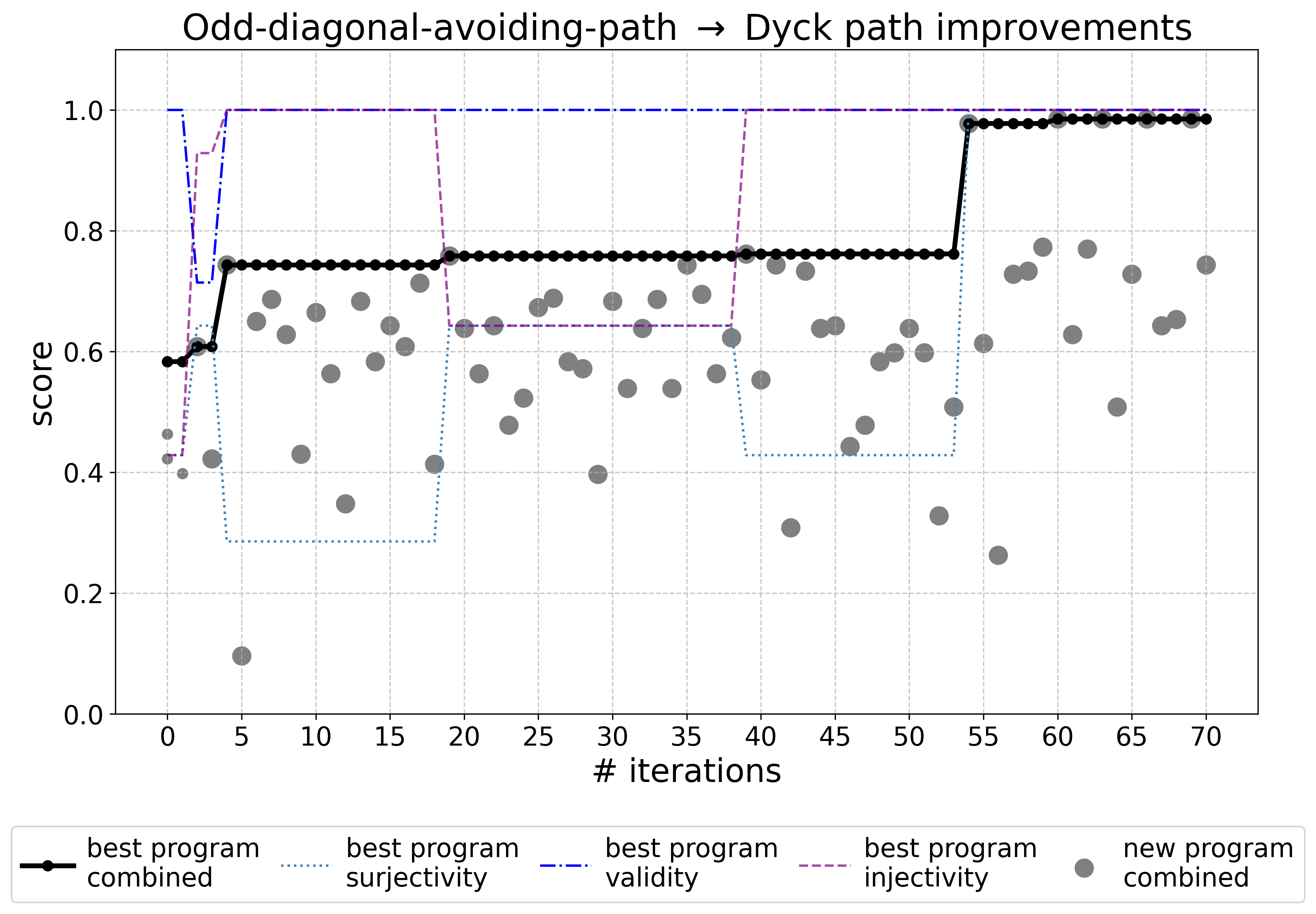}
\includegraphics[width=.925\textwidth]{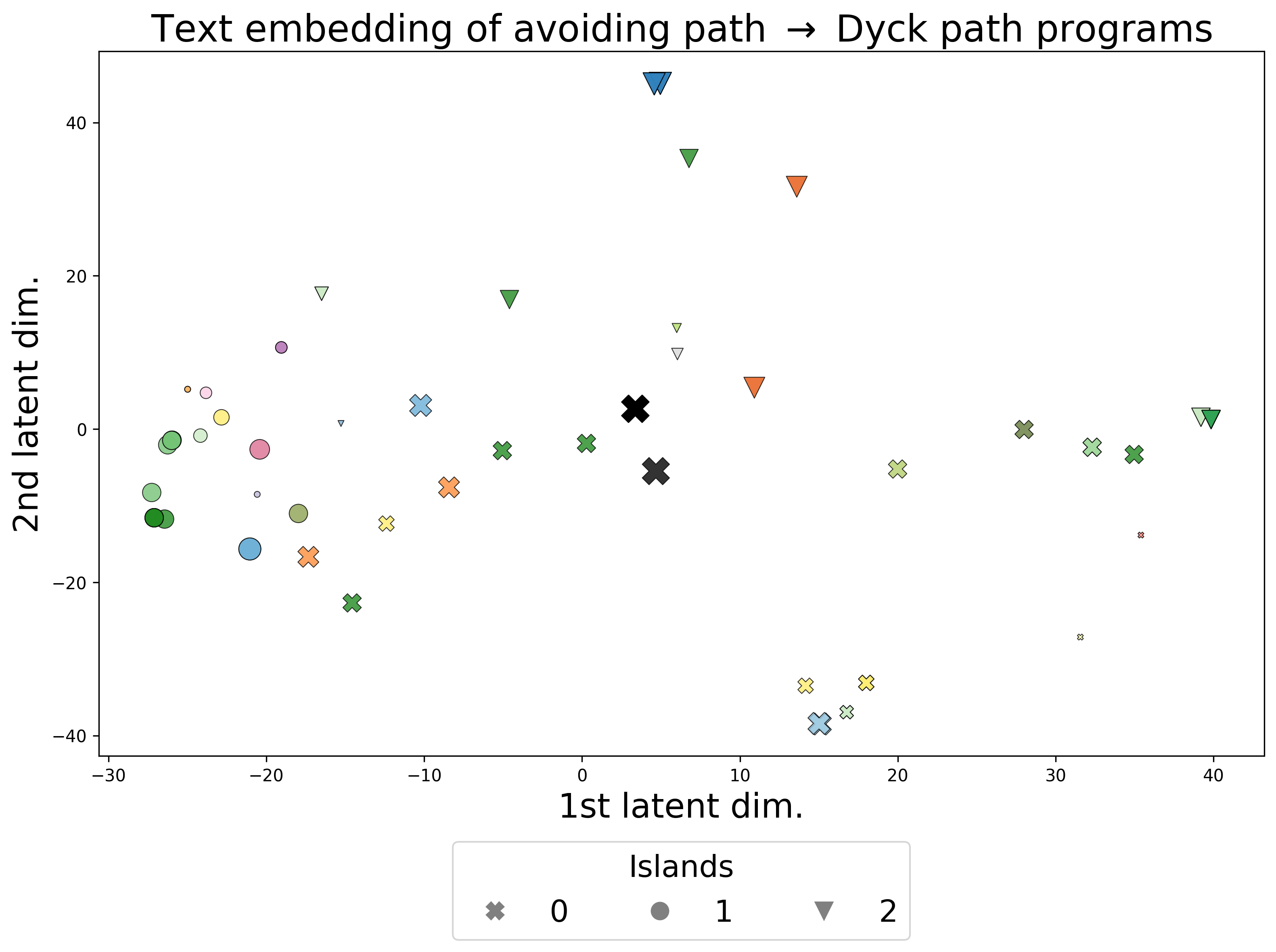}
\end{minipage}
\begin{minipage}{.475\textwidth}
\begin{pythonbox}
\begin{lstlisting}[style=pythonstyle]
def avoiding_path_to_dyck_path(path, m):
  # Find the first return to the diagonal: smallest t>0 with height h(t) = 0.
  h = 0
  t = -1
  for i, step in enumerate(path, start=1):
    h += 1 if step == 1 else -1
    if h == 0:
      t = i
      break
  if t == -1:
    raise ValueError("No return to diagonal found; invalid input.")

  U = path[:t]
  V = path[t:]

  # Determine the Dyck interior X of the primitive excursion U.
  if U[0] == 1:
    X = U[1:-1]  # interior is a Dyck path
    R = avoiding_path_to_dyck_path(V, m-t //2)
    return [1] + X + [0] + R
  else:
    # U = 0 Y 1 with Y a co-Dyck interior
    # take X = complement(Y) which is Dyck
    Y = U[1:-1]
    X = [1 - s for s in Y]
    R = avoiding_path_to_dyck_path(V, m-t //2)
    return [1] + R + [0] + X
    \end{lstlisting}
    \end{pythonbox}
\end{minipage}
\caption{Left, top: Progression of the best program metrics over 70 iterations for different terms in the scoring function. Left, bottom: Text embedding of all programs, after applying dimensionality reduction using PCA. Each functionally distinct mapping is given a different color (the bijection is black). The surjectivity score is indicated by the relative size of the points. Right: The valid bijection produced after 60 iterations.}
\label{fig:avoiding_program_embedding}
\end{figure}

\paragraph{Results.}
An evolving population of programs generated by OpenEvolve (including a valid solution) is visualized in Figure~\ref{fig:avoiding_program_embedding}. The top left plot shows the progression of program metrics: the system discovers the bijection after 60 iterations, but the combined scores of new introduced programs are very noisy. The best program validity score within the population is maximized early in evolution, indicating that the LLMs are easily able to generate programs that reliably output Dyck paths. On the other hand, surjectivity and injectivity take many more generations to optimize, with some understandable trade-offs occurring between these. The program defining the valid bijection which is generated at iteration 60 can be found in Figure~\ref{fig:avoiding_program_embedding}, right. This is the same as the known bijection from \cite{stanley2015catalan}. Notably, the docstring was over 600 words and included a sketch of invertibility (see Figure~\ref{fig:avoiding_docstring} in the Appendix), illustrating how the docstring produced by the system could aid a human mathematician in proving a bijection, provided that it is carefully checked.

\paragraph{How much variation is there is proposed solutions?} Surprisingly, we found limited variation in injectivity, surjectivity, and validity scores suggesting that models tend to gravitate toward a handful of different candidate maps. Further analysis showed that these metrics were not granular enough to always differentiate distinct mappings -- the 73 programs produce 26 functionally different mappings but only 10 different surjectivity scores. To further understand the variability among programs with the same metrics and even the same empirical mappings, we plot the text embeddings\footnote{We map programs into $\mathbb{R}^{3072}$ using OpenAI's text-embedding-3-large model. Vector proximity reflects textual semantic similarity; programs with similar structure, variable names, and comments will be close.} of all 73 programs after dimensionality reduction (Figure~\ref{fig:avoiding_program_embedding}, bottom left), giving programs defining the same mapping the same color. We find that programs with the same surjectivity score and even those that functionally correspond to the same map could have very different text embeddings (for example, the two orange upside down triangles that are the exact same program are distant on the plot). This highlights an axis of variation that can be easy to overlook: a significant amount of variation introduced by LLMs can relate to program implementation rather than program functionality.

\subsection{321-avoiding permutations and Dyck paths}

\paragraph{Problem.} Produce a program that implements a bijective map from Dyck paths of semilength $n$ to 321-avoiding permutations in $S_n$. 

\begin{figure}[!hbt]
    \centering
\begin{minipage}{.5\textwidth}
\includegraphics[width=.95\textwidth]{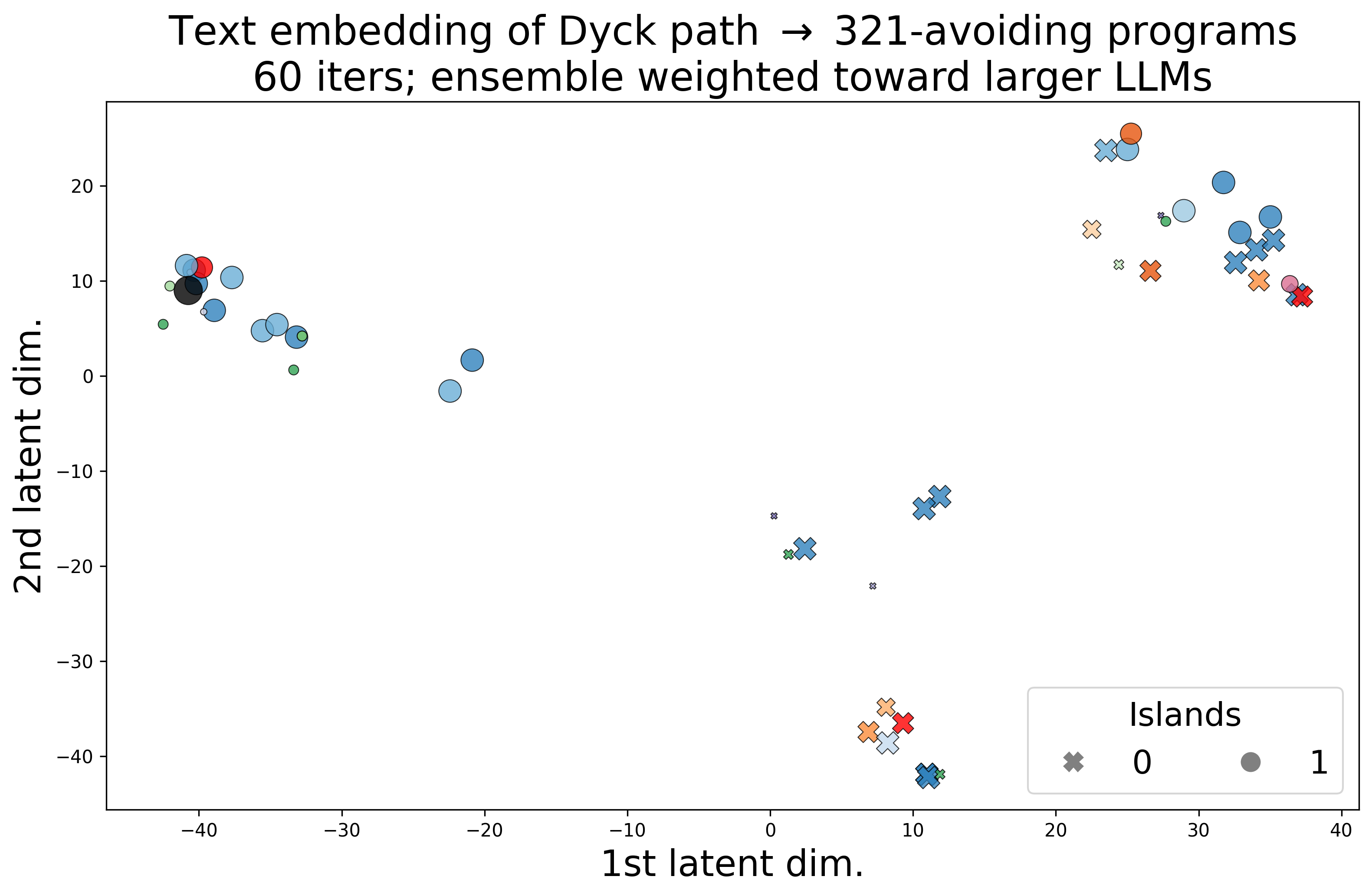}
\includegraphics[width=.95\textwidth]{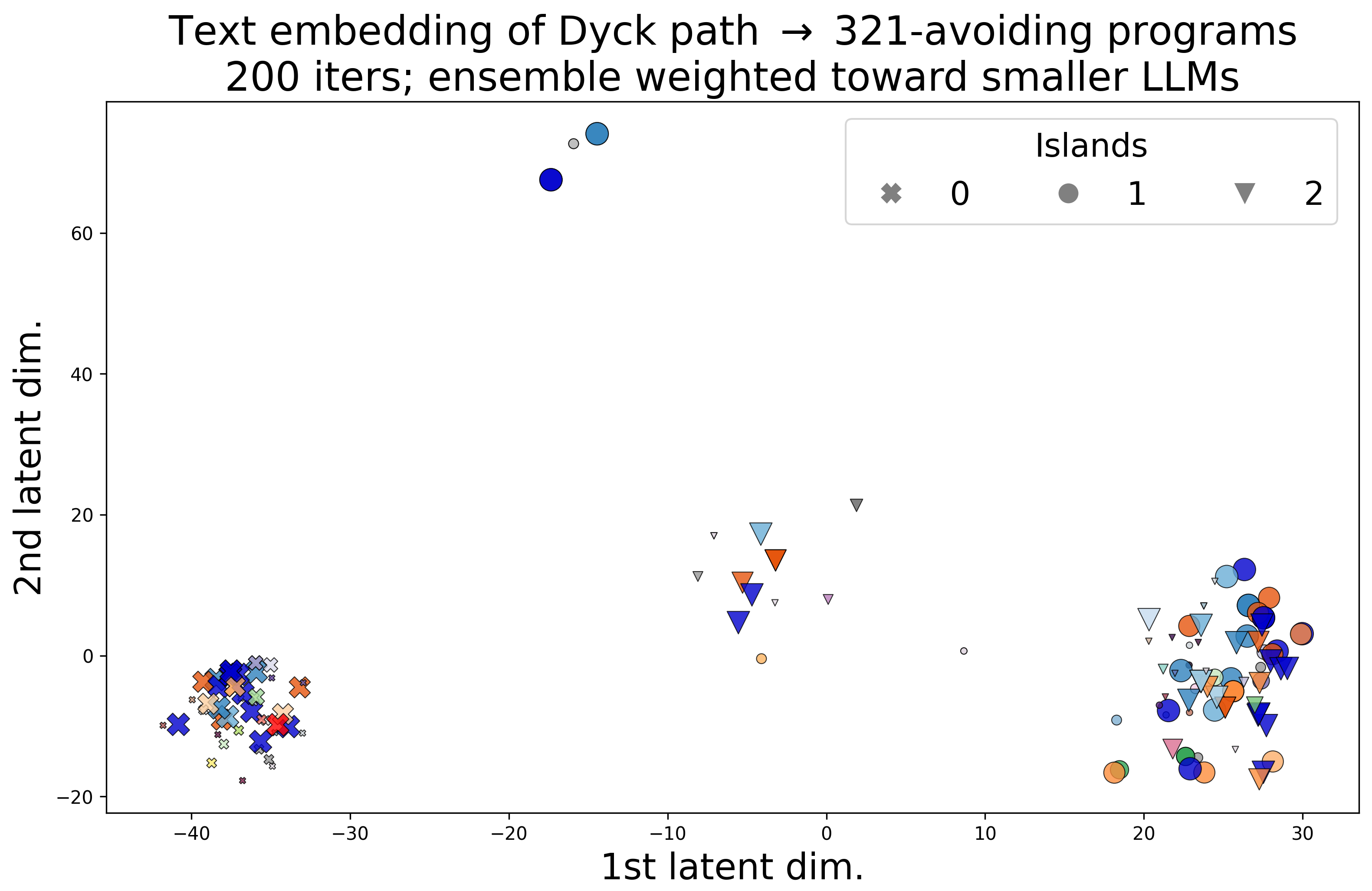}
\end{minipage}
\begin{minipage}{.45\textwidth}
\begin{pythonbox}
\begin{lstlisting}[style=pythonstyle]
def map(path):
  n = len(path) // 2
  if n == 0:
    return []
        
  perm = []
  S = []  
  current_num = 1 
  
  for step in path:
    if step == 1: 
      S.append(current_num)
      current_num += 1
    elif step == 0: 
      last = perm[-1] if perm else -1
      to_move = []
      while S and S[-1] < last:
        to_move.append(S.pop())
      if S:
        perm.append(S.pop())
      perm.extend(reversed(to_move))
      
  while S:
    perm.append(S.pop())
  return perm  
    \end{lstlisting}
    \end{pythonbox}
\end{minipage}
\caption{Left: Text embeddings of all programs, after dimensionality reduction using PCA (using more large models (\textbf{top}) or more smaller models (\textbf{bottom})). Functionally equivalent programs are given the same color. Surjectivity score is indicated by the relative size of the points. Right: The best program produced for the 321-avoiding permutation problem. The docstring and comments have been omitted for space.}
\label{fig:321_program_embedding}
\end{figure}

\paragraph{Results.} We applied OpenEvolve with a variety of configurations, but no run yielded a bijection. Figure~\ref{fig:321_program_embedding} shows results from two runs, one where the team of LLMs included more large models and one where the team included more small models\footnote{Smaller, cheaper models allow for more iterations (and therefore more candidate solutions) for the same cost.}.
In the former case, the best program mapped onto the subset of permutations that is both 321-avoiding and 3142-avoiding. For $n=4$, this program produces 13/14 target permutations, a significant improvement over the initial program which produces only 8/14 target permutations. In contrast, the team of smaller models failed to improve upon the empirical scores of the best initial program. 
In this case, we also tried evaluating on $n=5$ but found no improvement. Interestingly, in both runs some responses mention the Billey--Jockush--Stanley bijection \cite{billey1993some} by name, but attribute it to the wrong map.

\paragraph{LLMs can struggle when there are several closely related problems.} In addition to being in bijection with 321-avoiding permutations, Dyck paths are in bijection with stack-sortable (231-avoiding) permutations, and 132-, 213-, and 312-avoiding permutations. Bijections to 132-, 213-, and 312-avoiding permutations can be obtained from the bijection to stack-sortable permutations by simple operations (reversing the permutation, replacing each value $x$ in the permutation by $n+1-x$, or combining these operations, respectively). 
{\em The prominence of these known bijections appears to be a major barrier preventing the models from generating the bijection to 321-avoiding permutations.} The existence of these few prototypical solutions may partly explain the more densely clustered program text embeddings in Figure~\ref{fig:321_program_embedding} compared to Figure~\ref{fig:avoiding_program_embedding}. 

\subsection{Area-bounce exchanging bijection}

\paragraph{Problem.} Produce a program that implements a bijective map from Dyck paths to Dyck paths that exchanges the area and bounce statistics. 

\paragraph{Set-up details.} Since this is an open problem, we seeded the population with a map from Ayyer and Sundaravaradan \cite{ayyer2025area, qtcatalan-bijection}. The system is instructed that the map provided is a valid bijection that exchanges area and bounce on a subset of Dyck paths, and is asked to analyze its failures and propose a targeted modification to generalize it. Programs were scored using paths of semilength $4$. The Ayyer-Sundaravaradan mapping hits 8 of the 14 Dyck paths, is injective, and doesn't produce any out-of-domain outputs, since it returns None if the Dyck path input is not one of the paths the map works for.

\paragraph{LLMs can struggle when constructions in the problem look similar to constructions seen during pretraining.} 
Concepts that seem very salient in LLM's training data interfered with its ability to make progress on this problem. For example, the definition of bounce used in \cite{ayyer2025area} is formulated in a way that is slightly non-standard, so the initial feedback from the LLM evaluators was focused on re-defining the bounce statistic. We encountered a similar issue due to the prominence of the zeta map in related literature. The zeta map sends the statistics (area, dinv) to (bounce, area). While this does not solve the problem, the models frequently stated that it did and replaced the entire program with a program that used it. We addressed both of these issues with careful prompting. Ultimately, the model produced a higher scoring solution but not using an algorithmic bijection, which we discuss further in Section~\ref{sec:lessons}. 

\section{Lessons learned}
\label{sec:lessons}

\subsection*{Challenge \#1: Designing an objective function}
\label{sec:challenge1}

{\bf A primary challenge in applying evolutionary search for bijection discovery is objective function design.} Combinatorial bijections are all-or-nothing constructions: either a map is a bijection or it is not. As such, the community has spent much less time thinking about how to quantify the extent to which an arbitrary map is a bijection. The ability to construct granular scoring functions (ideally continuous) was identified as a crucial ingredient to finding solutions with AlphaEvolve in \cite{georgiev2025mathematical}.
In our case studies, programs that defined distinct mappings yielded identical surjectivity and injectivity scores, resulting in a plateaued objective landscape. 

It may be that the problem of bijection discovery is simply less amenable to incremental improvement than problems like sphere packing, as the discovery of a correct mapping can require a single critical insight that resolves multiple structural flaws at once. Alternatively, program evaluation may simply need refinement with more granular objective functions, such as a score that differentiates between a mapping with many 2-way collisions and one with a single $k$-way collision. It's also possible that a more fundamental change may be necessary. In \cite{georgiev2025mathematical}, they evolved programs that search for a construction rather than programs that directly generate a construction. Perhaps performance in our context could similarly be improved by evolving programs whose purpose is to search for bijections rather than evolving programs that define bijections.

\subsection*{Challenge \#2: Aligning the objective with mathematical intent}

{\bf Models often discover ``loopholes'' in the objective function, producing programs that achieve high scores but are not useful from the perspective of mathematical discovery.} This is referred to as the ``cheating phenomenon'' \cite{georgiev2025mathematical} or ``reward hacking'' \cite{skalse2022defining, lange2025shinkaevolve}, and occurs when the solution satisfies the literal constraints but violates the spirit of the problem. As we discussed in Section~\ref{sec:background}, one cheating behavior that we observed was LLMs exhaustively enumerating all objects from the sets $A$ and $B$ in order to construct a trivial index-based mapping. 
To prevent this kind of solution, we added explicit instructions not to generate all of the objects in $A$ and $B$. Additionally, we instructed the LLM evaluator to check for cheating, giving specific examples of what this might look like.

While these steps successfully prevented this particular instance of reward hacking, we found that it was difficult to anticipate all of the different ways a model might reward-hack. For instance, in evolving the area-bounce exchanging bijection (Problem 3), the system implemented a function that, given a Dyck path, first found its area and then used breadth-first search to find a path with the same bounce. While this program resulted in a near-perfect evaluation score, a mathematician would immediately disqualify such a solution, since it is a search algorithm rather than an algorithmic bijection and fails to provide any structural insight.  

For the foreseeable future, we believe that in bijection discovery problems an expert-in-the-loop is essential. Their role is not just to define the initial problem and objective, but also to periodically inspect the top-performing solutions for pathological behavior and evaluate whether the system is capturing the intention of the problem.  

\subsection*{Challenge \#3: The limitations of LLMs as evaluators}

The success of these systems currently depends on having objective, machine-verifiable evaluation metrics. While LLMs can be used to evaluate programs, two issues prevent them from being reliable evaluators in this context: the difficulty of capturing mathematical nuances in prompts, and the inconsistency of their responses. 

Aside from the case studies described in this paper, we also attempted to apply OpenEvolve to the problem of discovering a {\em combinatorial interpretation}: given a formula, produce a program that constructs a family of combinatorial objects counted by this formula. 
While we could automatically verify the correctness of the counts of candidate programs, we were not able to get reliable LLM evaluations of the quality of the interpretations. 
This was likely because (1) it is difficult to define what it means to be an interesting and useful combinatorial interpretation, making the evaluation prompt difficult to write, and (2) the candidate programs are not easy to evaluate even for a human. 
The candidate programs often contained long docstrings with elaborate definitions that seemed promising, but were ill-defined or nonsensical once one tried to actually construct an example of the object. This highlights a point that is perhaps obvious but still worth stating: {\bf given that these systems generate hundreds of programs, it is crucial to set up the problem and evaluation in such a way that doesn't require manual reading of LLM responses.} 

Beyond the challenge of defining subjective criteria, we also observed that LLM-based scoring was not consistent. 
While it might be expected that numerical evaluations vary depending on the model, we were surprised to find that even when the model was fixed, repeated evaluations of the exact same program gave scores spanning the entire possible range. This unreliability made it challenging to weight the LLM-based portion of the evaluation. When the LLM evaluation is weighted too low, pathological solutions are not sufficiently penalized, but when it is weighted too high, scores that are spuriously high due to inherent randomness interfere with the search. Given these issues, we think a better use of LLMs in the evaluation pipeline would be as filters (for instance, to automatically remove ``cheating'' programs) rather than as scorers.

\subsection*{Challenge \#4: LLMs favor solutions to well-known problems}

Beyond the observation that evolutionary program synthesis works best on problems with a smooth objective landscape, a more subtle challenge arises from biases in the LLM's training data. Our case studies illustrate that models have a tendency to get stuck on famous mathematical results that are associated with key words in the problem, even when they repeatedly prove to be unhelpful. For instance, when tasked with discovering a bijection between Dyck paths and 321-avoiding permutations (Problem 2) the model repeatedly produced programs that involved stack-sorting, which is used in well-known bijections with Dyck paths for four other pattern-avoiding permutations. Similarly, in Problem 3, the model was drawn to using the zeta map. These observations suggest that in problem selection and prompt writing, one should consider potential ``gravitational pulls'' and make efforts to steer models away from unproductive paths.

\subsection*{Challenge \#5: Overhead of setup}

Practical barriers to entry include the robustness and usability of the available software, and the overhead of problem formulation. The technical barrier is gradually lowering as codebases like OpenEvolve are regularly updated and new options are released. As the available software improves, the primary overhead will likely be problem selection, designing the objective function, and prompt engineering. The latter might be the least familiar for mathematicians new to using these tools, and we found that generating good prompts required significant iteration at first. In general, our experience was that getting these systems set up with a new problem requires at least a few hours even when you are familiar with the system and the problem, and longer otherwise.

{\bf Cost.} Finally, since these systems work best with frontier models requiring paid API calls, cost is an obvious concern. Fortunately, cost-per-token at time of writing means that significant compute is within reach of most mathematicians. For example, one run for Problem 2 used Gemini-2.5-flash, o4-mini, Gemini-2.5-pro, and Claude-3.7-sonnet for evolution and Gemini-2.5-flash for evaluation. A run of 200 iterations cost about \$10, with each iteration using two API calls: one for evolution (11,000 tokens on average) and one for evaluation (6,000 tokens on average).  \\

In summary, we remain optimistic about the application of evolutionary systems for bijection discovery, but successful application will require the right kind of problem with a well-designed objective function along with human expertise and guidance. One interesting additional consequence of working with these systems is that they force a mathematician to very directly confront what they value in a mathematical artifact. Monitoring LLMs that have an endless ability to produce mathematics that aligns with our literal instructions but not with what we want has the potential to sharpen our focus on the value of what we produce in our profession.

\section{Acknowledgements and Disclosure of Funding}
This work was conducted under the Laboratory Directed Research and Development Program at
PNNL, a multi-program national laboratory operated by Battelle for the U.S. Department of Energy
under contract DE-AC05-76RL01830.




\bibliography{sample}

@article{ayyer2025area,
  title={An area-bounce exchanging bijection on a large subset of {Dyck} paths},
  author={Ayyer, Arvind and Sundaravaradan, Naren},
  journal={Annals of Combinatorics},
  pages={1--29},
  year={2025},
  publisher={Springer}
}

@misc{cheng2025barbarians,
  title={Barbarians at the Gate: How {AI} is Upending Systems Research},
  author={Cheng, Audrey and Liu, Shu and Pan, Melissa and Li, Zhifei and Wang, Bowen and and Krentsel, Alex and Xia, Tian and Cemri, Mert and Park, Jongseok and Yang, Shuo and others},
 eprint={2510.06189},
 eprinttype={arxiv},
  year={2025}
}

@misc{liu2024llm4ad,
  title={{LLM4AD}: A platform for algorithm design with large language model},
  author={Liu, Fei and others},
  eprint={2412.17287},
 eprinttype={arxiv},
  year={2024},
note={arXiv:2412.17287}
}

@article{romera2024mathematical,
  title={Mathematical discoveries from program search with large language models},
  author={Romera-Paredes, Bernardino and others},
  journal={Nature},
  volume={625},
  number={7995},
  pages={468--475},
  year={2024},
  publisher={Nature Publishing Group UK London}
}

@misc{ellenberg2025generative,
  title={Generative modeling for mathematical discovery},
  author={Ellenberg, Jordan S and Fraser-Taliente, Cristofero S and Harvey, Thomas R and Srivastava, Karan and Sutherland, Andrew V},
eprint={2503.11061},
 eprinttype={arxiv},
note={arXiv:2503.11061},
  year={2025}
}

@misc{lange2025shinkaevolve,
  title={Shinka{E}volve: Towards open-ended and sample-efficient program evolution},
  author={Lange, Robert Tjarko and Imajuku, Yuki and Cetin, Edoardo},
eprint={2509.19349},
 eprinttype={arxiv},
note={arXiv:2509.19349},
  year={2025}

}

@article{alphaevolve,
title={{AlphaEvolve}: A {Gemini-powered} coding agent for designing advanced algorithms},
author={Novikov, Alexander and Vu, Ngân and Eisenberger, Marvin and Dupont, Emilien and Huang, Po-Sen and Wagner, Adam Zsolt and Shirobokov, Sergey and Kozlovskii, Borislav and Ruiz, Francisco J. R. and Mehrabian, Abbas and Kumar, M. Pawan and See, Abigail and Chaudhuri, Swarat and Holland, George and Davies, Alex and Nowozin, Sebastian and Kohli, Pushmeet and Balog, Matej},
year={2025},
url={https://deepmind.google/discover/blog/alphaevolve-a-gemini-powered-coding-agent-for-designing-advanced-algorithms/}
}

@book{stanley2015catalan,
  title={Catalan numbers},
  author={Stanley, Richard P},
  year={2015},
  publisher={Cambridge University Press}
}

@article{grosz2023bijectionist,
  title={A bijectionist’s toolkit},
  author={Grosz, Alexander and Kietreiber, Tobias and Pfannerer, Stephan and Rubey, Martin},
  journal={S{\'e}minaire Lotharingien de Combinatoire},
  volume={89},
  year={2023}
}

@misc{findstat,
author        = {Martin Rubey and Christian Stump and others},
title         = {{FindStat} - {T}he combinatorial statistics database},
url = {http://www.FindStat.org}
}

@article{elder2024homomesies,
  title={Homomesies on permutations: An analysis of maps and statistics in the FindStat database},
  author={Elder, Jennifer and Lafreni{\`e}re, Nadia and McNicholas, Erin and Striker, Jessica and Welch, Amanda},
  journal={Mathematics of Computation},
  volume={93},
  number={346},
  pages={921--976},
  year={2024}
}

@article{dowling2023homomesy,
  title={Homomesy on permutations with toggling actions},
  author={Dowling, William and Lafreniere, Nadia},
  journal={Involve},
volume={18},
year={2025},
pages={829–854}
}

@misc{adams2024cyclic,
  title={Cyclic sieving on permutations: an analysis of maps and statistics in the {FindStat} database},
  author={Adams, Ashleigh and Elder, Jennifer and Lafreniere, Nadia and McNicholas, Erin and Striker, Jessica and Welch, Amanda},
 eprinttype={arxiv},
eprint={2402.16251},
  year={2024},
note={arXiv:2402.16251}

}

@misc{openevolve,
  author = {Sharma, Asankhaya},
  title = {{OpenEvolve}: an open-source evolutionary coding agent},
  howpublished = {GitHub repository},
  url= {https://github.com/codelion/openevolve},
  note = {Accessed: 2025-08-10}
}

@misc{qtcatalan-bijection,
  author = {Sundar, Naren and Ayyer, Arvind},
  title = {qtcatalan-bijection},
  howpublished = {GitHub repository},
  url = {https://github.com/nanonaren/qtcatalan-bijection
}
}

@misc{callan2007bijections,
  title={Bijections from {Dyck} paths to 321-avoiding permutations revisited},
  author={Callan, David},
eprint={0711.2684},
 eprinttype={arxiv},
  year={2007},
note={arXiv:0711.2684}
}

@article{billey1993some,
  title={Some combinatorial properties of {Schubert} polynomials},
  author={Billey, S.C. and Jockusch, William and Stanley, Richard P},
  journal={Journal of Algebraic Combinatorics},
  volume={2},
  number={4},
  pages={345--374},
  year={1993},
  publisher={Springer}
}

@manual{sage,
  Key          = {Sage},
  Author       = {W.\thinspace{}A. Stein and others},
  Organization = {The Sage Development Team},
  Title        = {{S}age {M}athematics {S}oftware ({V}ersion 10.6)},
  url = {http://www.sagemath.org}
}

@misc{georgiev2025mathematical,
  title={Mathematical exploration and discovery at scale},
  author={Georgiev, Bogdan and G{\'o}mez-Serrano, Javier and Tao, Terence and Wagner, Adam Zsolt},
eprint={2511.02864},
  eprinttype={arxiv},
  year={2025},
note={arXiv:2511.02864},
}

@misc{Mouret2015IlluminatingSS,
  title={Illuminating search spaces by mapping elites},
  author={Jean-Baptiste Mouret and Jeff Clune},
eprint={1504.04909},
  eprinttype={arxiv},
  year={2015},
note={arXiv:1504.04909}
}

@misc{oeis,
    Author = {{OEIS Foundation Inc.}},
    Note = {Published electronically at \url{http://oeis.org}},
    Title = {The {O}n-{L}ine {E}ncyclopedia of {I}nteger {S}equences}}

@misc{nikankin2024arithmetic,
  title={Arithmetic without algorithms: Language models solve math with a bag of heuristics},
  author={Nikankin, Yaniv and Reusch, Anja and Mueller, Aaron and Belinkov, Yonatan},
eprint={2410.21272},
  eprinttype={arxiv},
  year={2024},
note={arXiv:2410.21272}
}

@article{skalse2022defining,
  title={Defining and characterizing reward gaming},
  author={Skalse, Joar and Howe, Nikolaus and Krasheninnikov, Dmitrii and Krueger, David},
  journal={Advances in Neural Information Processing Systems},
  volume={35},
  pages={9460--9471},
  year={2022}
}
\bibliographystyle{abbrv}

\newpage
\appendix
\section*{Appendix}
\addcontentsline{toc}{section}{Appendix}

\section{The Billey-Jockush-Stanley bijection}
We describe the Billey-Jockush-Stanley bijection \cite{billey1993some}, following the exposition in \cite{callan2007bijections}. Given a Dyck path defined as a sequence of $N$ and $E$ steps, the {\em ascent sequence} $\{a_i\}$ is the sequence that gives the numbers of the consecutive $N$ steps. The {\em descent sequence} $\{d_i\}$ is similarly defined by giving the numbers of consecutive $E$ steps. Then, define the {\em ascent code} (resp. {\em descent code}) to be the sequence of partial sums of the ascent lengths $A_i := \sum\limits_{j=1}^{i} a_j$ (resp. descent lengths $D_i := \sum\limits_{j=1}^{i} d_j$). The last element of the ascent and descent codes are omitted since it is always the length of the path. In a permutation $\sigma$, an {\em excedance location} is a value $i$ such that $\sigma(i) > i$, and $\sigma(i)$ is called the {\em excedance value}. With these definitions established, the bijection is straightforward to describe: given a Dyck path, let $A$ and $D$ be its ascent and descent code, respectively. Let $D$ be the excedance locations, and let $A+1$ be the excedance values. Fill in the remaining entries of the permutation in increasing order.

\section{Odd-diagonal avoiding paths and Dyck paths}

\subsection{Impact of model choice} Our initial run used an LLM ensemble for evolution with Gemini-2.5-pro, GPT-5, and Gemini-2.5-flash. Since GPT-5 was the model that produced the bijection, we did another run using just GPT-5 and found that using GPT-5 alone was sufficient, even with low reasoning. We subsequently did several runs using just Gemini-2.5-pro and found that the model produced the bijection as well when using high reasoning.

\subsection{Docstring of best program}

\begin{figure}[htb]
\begin{markdownbox}
\begin{lstlisting}[style=markdownstyle]
Direct, local bijection by first even-return decomposition with signed interior.

Input/Output:
- Steps are encoded by 1 = North (N), 0 = East (E).
- The input is an odd-diagonal-avoiding (ODA) path from (0,0) to (m,m), i.e., it never visits any diagonal point (2i-1, 2i-1). Equivalently, if h(t) is the prefix height (#1 - #0). Since the path avoids all diagonal points (2i-1, 2i-1), any return to the diagonal must occur at an even-coordinate point (2k, 2k).This means the prefix path to that point has length 4k, so any return to zero height must occur at a step index t that is a multiple of 4.
- The output is a Dyck path of semilength m: a 0/1-list of length 2m whose prefix heights are never negative.

Bijection overview:
- Decompose the input path P at its first return to the diagonal (which occurs at an index t |\(\equiv\)| 0 (mod 4) because P is ODA):
  P = U |\(\cdot\)| V,  with U a primitive excursion (no diagonal touch inside).
There are exactly two types of primitive excursions U from the diagonal to itself:
  (A) above-diagonal:     U = 1 X 0, with X a Dyck path,
  (B) below-diagonal:     U = 0 Y 1, with Y a `co-Dyck' (i.e., complement of a Dyck path).
In case (B), the complemented interior X := comp(Y) (swap 0 |\(\leftrightarrow\)| 1 in Y) is a Dyck path.
Thus, in both cases we have canonically associated a Dyck path X to the interior of U.

Recursive construction (phi):
- Base: if m = 0 (empty path), return [].
- Find t = first index with height h(t) = 0.
  Let U = path[:t], V = path[t:].
  Let X be the Dyck interior of U defined as:
        if U[0] == 1 then X = U[1:-1]           (above-diagonal excursion)
        if U[0] == 0 then X = comp(U[1:-1])     (below-diagonal excursion)
      Then define
        if U[0] == 1: phi(P) = 1 |\(\cdot\)| X |\(\cdot\)| 0 |\(\cdot\)| phi(V)
        if U[0] == 0: phi(P) = 1 |\(\cdot\)| phi(V) |\(\cdot\)| 0 |\(\cdot\)| X
      This uses only local data at the first even return and recurses on the suffix V,
      which is again ODA (concatenation of ODA excursions).

Why the output is Dyck:
- The word 1 X 0 is a Dyck excursion because X is Dyck. Concatenating Dyck excursions at returns preserves the Dyck property, so both constructions yield Dyck paths.

Invertibility (sketch):
- Given a Dyck path D, write its unique first-return decomposition D = 1 A 0 B with Dyck A, B. There are two mutually exclusive cases coming from the forward map:
  (i) D = 1 X 0 phi(V)       (above case), where X = A and phi(V) = B,
  (ii) D = 1 phi(V) 0 X      (below case), where phi(V) = A and X = B.
The inverse identifies which case holds by recursively inverting the left (A) and right (B) Dyck components and checking which reconstruction produces an ODA primitive excursion first (length of that excursion must be a multiple of 4 and its interior must be, respectively, Dyck (above) or complemented Dyck (below)). This test is local and deterministic, giving a two-way inverse and hence a bijection.

Algorithmic steps implemented here:
1) Validate length and balance.
2) If m == 0: return [].
3) Scan once to find the first index t > 0 with prefix height h(t) = 0 (first even return).
4) Set U = path[:t], V = path[t:].
  - If U starts with 1: X = U[1:-1].
  - If U starts with 0: X = complement(U[1:-1]) (bitwise 1 - s).
5) Recursively compute R = phi(V).
6) Return [1] + X + [0] + R if U[0] == 1, else [1] + R + [0] + X.
\end{lstlisting}
\end{markdownbox} 
\caption{The docstring of the program in Figure~\ref{fig:avoiding_program_embedding}, describing the avoiding path to Dyck path bijection.}
\label{fig:avoiding_docstring}
\end{figure}



\end{document}